\documentclass[11pt]{amsart}

\usepackage{amsmath,amsthm,amsfonts,amssymb,bm,graphicx }

\usepackage
{hyperref}
\hypersetup{colorlinks=true,citecolor=blue,linkcolor=blue,urlcolor=blue,
pdfstartview=FitH, pdfauthor=Valentin Blomer and P\'eter Maga} 

\theoremstyle{plain}
\newtheorem{theorem}{Theorem}
\newtheorem{lemma}{Lemma}

\newtheorem{cor}[lemma]{Corollary}

\numberwithin{equation}{section}
\numberwithin{theorem}{section}
\numberwithin{lemma}{section}
 
\newcommand{\LL}{L}

\begin{document} 

 \author{Valentin Blomer}
 \author{P\'eter Maga}
\address{Mathematisches Institut, Bunsenstr. 3-5, 37073 G\"ottingen, Germany} \email{blomer@uni-math.gwdg.de}  \email{pmaga@uni-math.gwdg.de}
\date{}

\title{The sup-norm problem for ${\rm PGL}(4)$}



\begin{abstract} Let $F$ be an $L^2$-normalized Hecke-Maa{\ss} cusp form for the group $\Gamma = {\rm SL}_{4}(\Bbb{Z})$  with Laplace eigenvalue $\lambda_F$. Assume that $F$ satisfies the Ramanujan conjecture at $\infty$ (this is satisfied by almost all cusp forms). If $\Omega$ is a compact subset of $\Gamma\backslash \mathcal{H}_4$, we show the bound $\| F|_{\Omega}\|_{\infty} \ll \lambda_F^{3/2 - \delta}$ for some absolute constant $\delta > 0$. 
\end{abstract}


\maketitle

\section{Introduction}

An automorphic form $F$ lives on a quotient $\Gamma \backslash X$ of a Riemannian symmetric space by a discrete subgroup of its  isometries, and it is in particular an eigenfunction of the Laplacian.  One of the   fundamental properties of an automorphic  form is its size and the distribution of its mass. Motivated by the correspondence principle in quantum mechanics, one is interested in sequences of (``high-energy'') eigenfunctions with increasing Laplace eigenvalue, and a central question is whether these eigenfunctions become in a certain sense equidistributed  on the space $\Gamma \backslash X$. In particular, the Quantum Unique Ergodicity conjecture of  Rudnick and Sarnak \cite{RuSa} asks  whether the measure $|F(z)|^2\, d\mu(z)$ converges to the uniform measure $d\mu(z)$ on $\Gamma \backslash X$ if the Laplace eigenvalue $\lambda_F$ of $F$ tends to infinity.     For detailed surveys of recent results we refer to \cite{Sa, Ze}.

Another measure of equidistribution is a bound of   some $L^p$-norm of $F$ for $p > 2$, the strongest case being the case $p = \infty$, because good upper bounds for $\| F \|_{\infty}$  rule out the possibility of high concentration of mass. If $X$ is a \emph{compact} locally  symmetric space of rank $r$ and $F$ is a joint eigenfunction  of the full algebra of   differential operators that are invariant under the Riemannian isometry group of $X$ (a polynomial algebra of rank $r$ including the Laplace operator), then one has the ``generic'' bound \cite{Sa1}
\begin{equation}\label{generic}
\| F\|_{\infty} \ll \lambda_F^{(\dim X -r)/4}.
\end{equation}
The compactness assumption is crucial here, as Brumley and Templier \cite{BT} have pointed out recently: for $n$ sufficiently large ($n \geq 6$ suffices), the bound  \eqref{generic} becomes false for  $X = {\rm PGL}(n)/{\rm PO}(n)$, because joint eigenfunctions have high peaks close to the cusps due to a certain degenerate behaviour of special functions. This phenomenon is already visible in weak form  in the classical case of ${\rm GL}(2)$ \cite{Sa1}, and a non-archimedean version can be found in \cite{Te}. Nevertheless, \eqref{generic} remains true if $F$ is restricted to a compact domain  of $\Gamma \backslash X$ so that the cuspidal regions are excluded.

While \eqref{generic} is sharp is general, many classical examples of Riemannian locally symmetric spaces 
enjoy additional symmetries given by the Hecke operators, a commutative family of normal operators. It is reasonable to assume that in this case \eqref{generic} can be improved by some power of $\lambda_F$ which would constitute a quantitative version of an equidistribution result.  The archetypical result of this kind is due to Iwaniec and Sarnak \cite{IS}  in the classical situation $X =  \mathcal{H}_2$, the hyperbolic plane, and $\Gamma \leq {\rm SL}_2(\Bbb{R})$ is a compact arithmetic subgroup or ${\rm SL}_2(\Bbb{Z})$. For $L^2$-normalized Hecke-Maa{\ss} cusp forms $F$ they proved the bound $\| F \|_{\infty} \ll \lambda_F^{5/24 +\varepsilon}$. 

The sup-norm problem for arithmetic manifolds has much in common with the subconvexity problem for automorphic $L$-functions in analytic number theory, in terms of the nature of the problem, the results available and the methods involved. Both the subconvexity and the sup-norm problem ask for a power saving in a specific arithmetic situation where purely analytic tools provide a ``generic" bound. Both for the subvonvexity and the sup-norm problem the situation for ${\rm GL}(2)$ is essentially well-understood while higher rank results are extremely rare and sporadic. For the sup-norm problem on ${\rm GL}(2)$ we mention in particular the works \cite{HT} for the level aspect,  \cite{Te1} for a hybrid version, \cite{Sah} for non-squarefree levels, \cite{Ki} for metaplectic forms and \cite{BHM} for forms over an imaginary quadratic number field. The first example for real rank  $> 1$ was the case of the symplectic group ${\rm Sp}(4)$ by the first author and A.\ Pohl \cite{BP}.  The case of ${\rm PGL}(3)$ was worked out simultaneously to the present paper in \cite{HRR}. (It may be remarked that  all of the results mentioned in this paragraph are very recent which shows the enormous dynamic that the sup-norm problem has gained in the past  months and years.) 
Finally the methods involved in the subconvexity and the sup-norm problem are related in the sense that they typically rely on harmonic analysis and trace formulae on the one hand, and some arithmetic input like estimation of character sums or diophantine analysis on the other. In particular the sup-norm problem features a beautiful mixture of spectral analysis (behaviour of spherical functions), combinatorial analysis (Hecke algebras and $p$-adic groups) and diophantine analysis (counting integral matrices with prescribed diophantine properties). 

We will explain these three ingredients in more detail in the next section and proceed with the statement of our main result. Going beyond \cite{HRR},  we will solve the sup-norm problem for ${\rm PGL}(4)$. In comparison, current technology has not yet produced any subconvexity results for genuine $L$-functions for ${\rm GL}(n)$ with $n > 3$, so  we believe that the present article advances substantially   our understanding of  automorphic forms on higher rank groups.  In fact, except for the last chapter, the complete article is written for general $n$.  We highlight in particular Sections \ref{sec4} and \ref{gen} where an explicit amplification scheme for ${\rm GL}(n)$ is developed that may be useful in other situations. 

For details and definitions in the following theorem we refer to Section \ref{sec3}.

\begin{theorem}\label{thm1}   Let $F$ be an $L^2$-normalized Hecke-Maa{\ss} cusp form for the group ${\rm SL}_4(\Bbb{Z})$ with Laplacian eigenvalue $\lambda_F$. Assume that $F$ satisfies the  Ramanujan conjecture at $\infty$, i.e.\ its archimedean Langlands parameters are real. Let $\Omega$ be a fixed compact subset of ${\rm SL}_4(\Bbb{Z})\backslash \mathcal{H}_4$. Then
$$\| F|_{\Omega}\|_{\infty} \ll_{\Omega} \lambda_F^{\frac{3}{2} - \delta}$$
for some absolute constant $\delta > 0$.
\end{theorem}

As mentioned before, the restriction to a compact subset $\Omega$ is quite natural in view of the work of Brumley-Templier \cite{BT}. 
The condition on temperedness at infinity is of technical nature, because our current test function is not necessarily bounded away from zero at the non-tempered spectrum. It is believed that this condition is automatically satisfied, and it is known \cite{LM} that the set of Hecke-Maa{\ss} cusp forms violating the Ramanujan conjecture at $\infty$ has density zero in the set of all Hecke-Maa{\ss} cusp forms (when ordered by Laplacian eigenvalue). The constant $\delta$ is absolute, and one could easily produce a (quite reasonable) numerical value. We have made no attempt to optimize the constant, and in fact, the purpose of this article to use techniques that have the potential to generalize to higher rank situations.  The current record  is $\delta < 1/24$ for $n=2$ \cite{IS} and $\delta < 1/124$ for $n=3$ \cite{HRR}.  Our proof gives in fact the slightly stronger bound
$$\| F|_{\Omega}\|_{\infty} \ll_{\Omega} \prod_{1 \leq j < k \leq 4} (1 + |\mu_j - \mu_k|)^{1/2 - \delta}$$
where $(\mu_1, \ldots, \mu_4)$ are the archimedean spectral parameters of $F$.

\section{Methods and ingredients}\label{sec2}

We give here an informal description of the methods and   ingredients involved in the proof of Theorem \ref{thm1}. The general strategy follows the ground-breaking ideas of \cite{IS}. Instead of estimating $F(z)$  for $z \in \Omega$ individually, we consider a weighted spectral sum  
\begin{equation}\label{sum}
  \sum_{\varpi} A(\varpi) |F_{\varpi}(z)|^2 
  \end{equation}
where the sum runs over the constituents $\varpi$ of $L^2(\Gamma \backslash \mathcal{H}_n)$ (including Eisenstein series, so that the sum is in reality a combination of sums and integrals) and $A(\varpi)$ is a non-negative weight function with $A(\varpi_0) = 1$ for the specific automorphic representation $\varpi_0$ generated by a Hecke-Maa{\ss} cusp form $F = F_{\varpi_0}$ whose sup-norm we want to bound. Applying a pre-trace formula to the expression \eqref{sum}, we obtain a weighted sum over matrices that leads to a diophantine problem. 
 
The key to success is a clever construction of the function $A(\varpi)$. Let $G = {\rm PGL}(n) = NAK$ be the Iwasawa decomposition and $W$ the Weyl group. We can assume that each $\varpi$ is an eigenform of the algebra of differential operators on $\Gamma \backslash \mathcal{H}_n = \Gamma \backslash G/K$ and the algebra of Hecke operators at all primes $p$. Hence each $F_{\varpi}$ comes with parameters $\mu = \mu_{\varpi} = (\mu_1, \ldots, \mu_n) \in \mathfrak{a}^{\ast}_{\Bbb{C}}/W$ at infinity and $\alpha_p = \alpha_p(\varpi) = (\alpha_{p, 1}, \ldots, \alpha_{p, n})$ at each finite place $p$. Our test function will factor into an infinite part  and a finite part. At the infinite place, we choose a bump function $\tilde{f}$ on $\mathfrak{a}^{\ast}_{\Bbb{C}}/W $  with $\tilde{f}(\mu_{\varpi_0}) = 1$ and essential support in a ball of radius 1 about $\mu_{\varpi_0}$, i.e. $\tilde{f}(\mu)$ is rapidly decaying if $|\mu - \mu_{\varpi_0}| \geq 1$. Of course it would be desirable to have an even smaller support of $\tilde{f}$, but then the geometric side of the trace formula gets out of control. The geometric side features the inverse spherical transform $f$ of $\tilde{f}$, constructed as an integral of $\tilde{f}$ against the elementary spherical functions. Using the Paley-Wiener theorem for $G$  \cite{Ga} and some asymptotic analysis of elementary spherical functions \cite[Theorem 2]{BP} one obtains a rapidly decaying function $f : K \backslash G / K \rightarrow \Bbb{C}$ on the geometric side with $f(0) \ll \|\mu_{\varpi_0} \|^{n(n-1)/4}$. As $\| \mu_{\varpi_0} \|^2 \asymp \lambda_F$, this is the order of magnitude of the generic bound. This part of the argument works for very general groups $G$. 

The finite part of $A(\varpi)$ has a similar purpose. For a prime $p$ we would like to choose a function on the set of Satake parameters that localizes at  $\alpha_p(\varpi_0)$ and is non-negative everywhere. The classical choice of Iwaniec-Sarnak \cite{IS} in the case $n=2$ is
\begin{equation}\label{am}
A_{\varpi}(p ) :=  \bigl| \overline{\lambda_{\varpi_0}( p)}  \lambda_{\varpi}(p) +  \overline{\lambda_{\varpi_0}( p^2)}  \lambda_{\varpi}(p^2) \bigr|^2
 \end{equation}
where $\lambda_{\varpi}(n)$ are the usual (normalized) Hecke eigenvalues. 
This function (as a function of $\varpi$) is non-negative and bounded from below for $\varpi = \varpi_0$ by the Hecke relations. 
As in the archimedean case, the geometric side of the trace formula features the inverse spherical transform of \eqref{am}, a locally constant bi-$K_p$-invariant function, or expressed more simply,    a function on double cosets $\Gamma   \text{diag}(p^{a_1}, \ldots, p^{a_n}) \Gamma$, typically referred to as the Hecke algebra. 
This idea works in principle in general, but there are two problems. On the one hand, one needs to find suitable Hecke relations that ensure  $A_{\varpi_0}(p ) \gg 1$. This is not too hard, and has been shown in great generality in \cite{SV}. For $ {\rm PGL}(n)$ one can be very explicit and build a suitable relation using the Hecke operators (see Lemma \ref{lem1})
\begin{equation}\label{am1}
\Gamma   \text{diag}(p^{j}, 1, \ldots, 1) \Gamma, \quad 1 \leq j \leq n.
\end{equation}
We refer to Section \ref{sec4} for details and just mention that the linear relation will feature the Kostka numbers as coefficients that appear quite frequently in the theory of Schur polynomials.  The more serious problem is the multiplication of the square in \eqref{am}. In general the multiplication of two double cosets in the Hecke algebra is a complicated combinatorial problem. Luckily, in the case of \eqref{am1}, i.e.
$$\Gamma   \text{diag}(p^{j}, 1, \ldots, 1) \Gamma \cdot \Gamma   \text{diag}(p^{j},  \ldots, p^j, 1) \Gamma$$
where the second factor is the adjoint operator, this is not too hard to achieve (see Lemma \ref{lem2}). The final amplifier is then a \emph{sum} over various $A_{\varpi}( p)$: when \eqref{am} is averaged over $p$ \emph{inside} the absolute values, the amplifier is small for $\varpi \not= \varpi_0$.  As is often the case in number theory, the key to success is to play off  many places of the underlying number field against each other. 

The geometric side of the trace formula features now a diophantine problem, and after spectral and combinatorial analysis, this is the third important ingredient, and really the heart of the matter. In a nutshell, one has to count matrices $\gamma \in \text{Mat}(n, \Bbb{Z})$ satisfying 
\begin{equation}\label{co}
\gamma^{\top} Q \gamma = (\det\gamma)^{2/n} Q + \text{ very small error }
\end{equation}
where $Q \in \text{Mat}(n, \Bbb{R})$ is a fixed positive definite matrix depending on the point $z \in \mathcal{H}_n$ at which we want to bound $F_{\varpi_0}$. Here 
\begin{equation*}\label{co1}
\det \gamma = p_1^j p_2^{j(n-1)}
\end{equation*}
for $1 \leq j \leq n$ and primes $p_1, p_2 \asymp L$ of the same order of magnitude. We know in addition that 
\begin{equation}\label{co2}
\text{the second determinantal divisor, i.e.\ the greatest common divisor of all 2-by-2 minors, equals $p_2^j$.}
\end{equation} 
In small rank, \eqref{co2} can be ignored, but in higher rank, it becomes an essential feature of a very interesting and challenging counting problem. The benchmark one has to beat is $L^{j(n-1)}$ for the number of such matrices.  

Let $x_1, \ldots,  x_n \in \Bbb{Z}^n$ denote the columns of $\gamma$. 
Constructing the columns one at a time, the $j$-th column satisfies a quadratic condition as well as $j-1$ linear conditions with respect to the preceding (already fixed) columns. The condition \eqref{co} alone gives us therefore 
\begin{equation}\label{l1}
L^{j\{(n-2)  + (n-3) + \ldots + 1 + \varepsilon\}}
\end{equation}
 solutions. This is a reflection of the fact that  the compact group ${\rm SO}(n)$ satisfies ${\rm SO}(n)/{\rm SO}(n-1) \cong S^{n-1}$, so inductively one is looking at quadratic problems with fewer and fewer variables in each step. 
 This reasoning is made precise in Corollary \ref{count} below and suffices for a quick proof of Theorem \ref{thm1} for   $n\leq 3$  without much effort (for some non-optimized value of $\delta$), since in this case \eqref{l1} is strictly smaller (by a power of $L$) than $L^{j(n-1)}$.  For $n = 4$ we hit exactly the benchmark, so we need to exploit the additional   condition \eqref{co2} which means essentially that all columns are multiples of each other modulo $p_2^j$. We will solve the corresponding diophantine  problem   in Section \ref{4}.

\section{General Set-up}\label{sec3}

For the rest of this paper let $$G = {\rm PGL}_n(\Bbb{R}),  \quad K = {\rm PO}(n),  \quad \Gamma = {\rm SL}_n(\Bbb{Z}),  \quad \mathcal{H}_n \cong G/K.$$ The latter is the generalized upper half plane (see \cite{Go}), a connected manifold of dimension $(n-1)(n+2)/2$. Let $W \cong S_n$ be the Weyl group, $A$   the diagonal torus in $G$, and  $\mathfrak{a}$  the corresponding Lie algebra. It has a root system of type $A_{n-1}$, and we have $$\rho = \frac{1}{2} \sum_{j=1}^n (n+1 - 2j)e_j \in \mathfrak{a}^{\ast}$$
for the half-sum of positive roots where $e_j(\text{diag}(a_1, \ldots, a_n)) = a_j$.  
As usual, we denote by $C_{\rho}$ the convex hull of the points $\{w \rho \mid w \in W\}$. We have a spectral decomposition
\begin{equation}\label{spec}
L^2(\Gamma \backslash \mathcal{H}_n) = \int V_{\varpi} d\varpi = L^2_{\text{cusp}}(\Gamma \backslash \mathcal{H}_n) \oplus L^2_{\text{Eis}}(\Gamma \backslash \mathcal{H}_n),
\end{equation}
where each $V_{\varpi}$ is a one-dimensional space generated by an  eigenform $F_{\varpi}$ (potentially an Eisenstein series or an iterated residue thereof) of the Hecke algebra and the algebra of invariant differential operators. We defer a detailed discussion of the Hecke algebra to the next section.  At the archimedean place,  each $F_{\varpi}$ comes with $n$ spectral parameters $\mu = (\mu_1, \ldots, \mu_n)$ satisfying 
\begin{equation}\label{Lambda}
  \sum_{j=1}^n \mu_j = 0, \quad \{\mu_1, \ldots, \mu_n\} = \{\bar{\mu}_1, \ldots, \bar{\mu}_n\}, \quad \mu \in \mathfrak{a}^{\ast} + i C_{\rho}.
  \end{equation}
(In this normalization, the tempered spectrum is real.) In particular, the spectral parameters are  real or come in complex conjugate pairs. The Jacquet-Shalika bound \cite{JS}   implies for cusp forms more strongly that $|\Im \mu_j| \leq 1/2$, but we do not need this.  We denote the set of possible spectral parameters, i.e.\ the subset of $\mathfrak{a}_{\Bbb{C}}^{\ast}$ satisfying \eqref{Lambda}, by $\Lambda$.  

The Laplace eigenvalue of $F = F_{\varpi}$ is given by
\begin{equation}\label{laplace}
\lambda_{F} = \frac{1}{2}(\| \rho \|^2 + \| \mu \|^2) = \frac{n^3 - n}{24} + \frac{1}{2}(\mu_1^2 + \ldots + \mu_n^2).
\end{equation}
The spectral density at $\lambda = \sum_j \lambda_j e_j \in \mathfrak{a}^{\ast}$ is given by  
\begin{equation}\label{density}
\frac{1}{|\textbf{c}(\lambda)|^2} \ll  \frac{1}{|\tilde{\textbf{c}}(\lambda)|^2} :=  \prod_{1 \leq j < k \leq n} (1 + |\lambda_j - \lambda_k|) \ll 1 + \|\lambda \|^{n(n-1)/2}
\end{equation}
where  $\textbf{c}(\lambda)$ denotes the Harish-Chandra $\textbf{c}$-function (see \cite[Section 3]{LM}).  \\

As in \cite{BP} we choose the following archimedean test function. Fix some large $\mu  \in \mathfrak{a}^{\ast} $. Let $\psi$ be a fixed Paley-Wiener function on $\mathfrak{a}_{\Bbb{C}}^{\ast}$ such that $\psi$ is non-negative on $\mathfrak{a}^{\ast}$ and $\psi(0) = 1 $. Then we choose
\begin{equation}\label{deff}
\tilde{f}_{\mu}(\lambda) := \Bigl( \sum_{w \in W} \psi(\mu - w \cdot \lambda)\Bigr)^2.
\end{equation}
This is again a Paley-Wiener function that is non-negative on all possible spectral parameters $\lambda \in \Lambda$ (because $\psi(\bar{\lambda}) = \overline{\psi(\lambda)}$, hence the term inside the parenthesis is real), and obviously   
\begin{equation}\label{geq1}
\tilde{f}_{\mu}(\mu)     \geq 1
\end{equation} 
as well as 
\begin{equation}\label{decay}
  \tilde{f}_{\mu}(\lambda) \ll_A \max_{w \in W}(1 + \| \mu - w \cdot \lambda\|)^{-A}
\end{equation}
for $\lambda \in \mathfrak{a}^{\ast}$ and any $A \geq 0$. 
  By the Paley-Wiener theorem for $G$ \cite{Ga}, its spherical inversion $f_{\mu}$ is a compactly supported, smooth, bi-$K$-invariant function. 

Assume that we want to bound the sup-norm of a Hecke-Maa{\ss} form $F_{\varpi_0}$ occurring in \eqref{spec} with large spectral parameter $\mu \in \mathfrak{a}^{\ast}$ (i.e. satisfying the Ramanujan conjecture at $\infty$). Then we use this $\mu$ in the definition \eqref{deff}. The pretrace formula states
\begin{equation}\label{pt}
\int \tilde{f}_{\mu}(\mu_{\varpi})F_{\varpi}(x) \overline{F_{\varpi}(y)} d\varpi = \sum_{\gamma \in \Gamma} f_{\mu}(x^{-1} \gamma y), \quad x, y \in G,
\end{equation}
where we interpret $F_{\varpi}$ as $K$-invariant functions on $\Gamma \backslash G$. The left hand side contains in particular an $L^2$-normalized version of our preferred Hecke-Maa{\ss} cusp form $F_{\varpi_0}$ for which we want to prove Theorem \ref{thm1}. 

By the Harish-Chandra inversion formula \cite[Ch.\ IV]{He}
$$f(g) = \frac{1}{|W|} \int_{\mathfrak{a}^{\ast}} \tilde{f}(\lambda) \phi_{\lambda}(g) \frac{d\lambda}{|\textbf{c}(\lambda)|^2},$$
in connection with the bound  \cite[Theorem 2]{BP} for the elementary spherical function $\phi_{\lambda}$ and \eqref{decay}, we have the important upper bound 
\begin{equation}\label{spherical}
f_{\mu}(g) \ll \frac{1}{|\tilde{\textbf{c}}(\mu)|^2}\bigl(1 + \| \mu \| \| C(g)\|\bigr)^{-1/2}
\end{equation}
where $C : G \rightarrow \mathfrak{a}/W$ is the Cartan projection, so that 
\begin{equation}\label{Cartan}
  g = k_1 \exp(C(g)) k_2.
\end{equation}  

\section{Hecke operators }\label{sec4}

 In the following we equip $\Bbb{Z}^n$ with lexicographic order,  i.e.\ for $\textbf{a} = (a_1, \ldots, a_n)$, $\textbf{b} = (b_1, \ldots, b_n) \in \Bbb{Z}^n$ we write $\textbf{a} \leq \textbf{b}$ if $\textbf{a} = \textbf{b}$ or if there is some index $1 \leq j \leq n$ such that $a_i=  b_i$ for $1 \leq i \leq j-1$ and $a_j < b_j$. We write $|\textbf{a}| = \sum_j a_j$.

Let   $p$ be a prime.  For a  double coset $\Gamma \text{diag}(p^{a_1}, \ldots, p^{a_n}) \Gamma = \cup_j \Gamma M_j$ with $\textbf{a} = (a_1, \ldots, a_n) \in \Bbb{Z}^n$, we define the   Hecke operator
\begin{equation}\label{defHeck}
T_{\textbf{a}}(p ) := T_{\Gamma \text{diag}(p^{a_1}, \ldots, p^{a_n}) \Gamma} : F \mapsto    \sum_j F(M_j \,\, \cdot)
\end{equation}
for any function $F : \Gamma \backslash \mathcal{H}_n \rightarrow \Bbb{C}$. 
The vector space generated by these operators is a commutative algebra. The product  of two such operators corresponds to the multiplication of two double cosets $\Gamma A \Gamma = \bigcup_j \Gamma A_j$,  $\Gamma B \Gamma = \bigcup_k \Gamma B_k$:
\begin{equation}\label{multi}
T_{\Gamma A \Gamma } \circ T_{\Gamma B \Gamma} \colon F \mapsto  \sum_{j, k}F( 
A_jB_k \,\cdot\ )
=  \sum_{D} \alpha_{D} T_{\Gamma  D \Gamma}
\end{equation}
where $D$ runs through a system of representatives of double cosets contained in $\Gamma A \Gamma B \Gamma$ and $\alpha_D$ is the number of pairs $(j, k)$ such that $\Gamma D = \Gamma A_jB_k$. 

We have
$$\overline{T}_{\textbf{a}}( p) =    T_{-\textbf{a}}( p).$$
The standard Hecke operators are given by
$$T(p^k) = \sum_{\substack{  |\textbf{a}| = k\\ a_1 \geq a_2 \geq \ldots \geq a_n \geq 0}} T_{\textbf{a}}( p).$$

Note that $T_{(a, \ldots, a)}( p ) = \text{id}$ for all $a \in \Bbb{Z}$, and also $T_{\textbf{a}}(p ) = T_{\sigma(\textbf{a})}(p )$ for all $\sigma \in S_n$. Hence we can assume without loss of generality that $ a_1 \geq a_2 \geq \ldots \geq a_n $ (and even $a_n = 0$). We write $(\Bbb{N}_0^n)^{\ast}$ for the set of all $n$-tuples $\textbf{a} = (a_1, \ldots, a_n)$ in non-increasing order $ a_1 \geq a_2 \geq \ldots \geq a_n\geq 0$, and  for $\textbf{a} = (a_1, \ldots, a_n) \in (\Bbb{N}_0^n)^{\ast}$ we write $$v(\textbf{a}) = \sum_{j=1}^n ja_j.$$ 

If $F$ is an eigenform of the full Hecke algebra, we denote by $\lambda_{\textbf{a}}(p, F)$  and $\lambda(p^k, F)$ the eigenvalue with respect to the operators $T_{\textbf{a}}( p)$  and $T(p^k)$, respectively. In this normalization, the Ramanujan conjecture states $|\lambda(p, F)| \leq n p^{(n-1)/2}$, and our $T(p )$ is $p^{(n-1)/2} $ times the corresponding operator $T_p$  in \cite{Go}. We refer to standard texts such as \cite{AZ} or \cite{Fr} for more details on the Hecke algebra for ${\rm GL}(n)$.   

The Satake map $\omega$ provides an isomorphism between the $p$-part of the Hecke algebra and symmetric functions in $\Bbb{Q}[x_1^{\pm}, \ldots, x_n^{\pm}]^{S_n}$. An explicit description of this map can be found in \cite[(1.7)]{An} which we reproduce here.  Let $\textbf{a}  \in (\Bbb{N}_0^n)^{\ast}$  and suppose that the multiset $\{a_1, \ldots, a_n\}$ contains $t$ distinct elements occurring with multiplicities $k_1, \ldots, k_t$, then the image of $T_{\textbf{a}}(p )$ is the polynomial 
\begin{equation}\label{and}
\omega(T_{\textbf{a}}(p ))(\textbf{x}) = \frac{(1 - 1/p)^n}{p^{v(\textbf{a})}} \prod_{i=1}^t \prod_{j=1}^{k_t} (1-1/p^j)^{-1} \sum_{\sigma \in S_n} \sigma\Bigl(\textbf{x}^{\textbf{a}} \prod_{1 \leq i < j \leq n} \frac{x_i - x_j/p}{x_i - x_j}\Bigr)
\end{equation}
where $\sigma $ acts on the $x$-variables. Viewing a Hecke operator as a locally constant bi-$K_p$-invariant function, the Satake map is, up to scaling,   just the   spherical transform: if ${\bm \alpha} = (\alpha_1, \ldots, \alpha_n)$ are the $p$-Satake parameters of an eigenform $F$, then 
\begin{equation*}\label{eig}
\lambda_{\textbf{a}}(p, F) =    \omega(T_{\textbf{a}}( p))(p^{(n+1)/2}{\bm \alpha}).  
\end{equation*}

In general, \eqref{and} is not easy to evaluate, but we have the following basic result:
\begin{lemma}\label{basic}
 For $\textbf{a} \in (\Bbb{N}^n_0)^{\ast}$ we have 
$$\omega(T_{\textbf{a}}(p ))(\textbf{x})  = p^{-v(\textbf{a})} \sum_{\substack{\textbf{b} \in  \Bbb{N}_0^n,   |\textbf{b}| = |\textbf{a}|\\  \sigma(\textbf{b}) \leq \textbf{a} \text{ for some } \sigma \in S_n}} c_p(\textbf{a}, \textbf{b}) \textbf{x}^{\textbf{b}}$$
with certain coefficients $c_p(\textbf{a}, \textbf{b}) \in \Bbb{Z}[1/p]$ satisfying $c_p(\textbf{a},  \textbf{b}) = c_p(\textbf{a}, \sigma(\textbf{b}))$ for all $\sigma \in S_n$, $c_p(\textbf{a}, \textbf{a}) = 1$ and $c_p(\textbf{a}, \textbf{b}) \ll 1$. 
\end{lemma}

Indeed, the  condition $|\textbf{b}| = |\textbf{a}|$ is clear from degree considerations, $c_p(\textbf{a}, \textbf{b}) \in \Bbb{Z}[1/p]$ and the $S_n$-symmetry  follow from the definition, $c_p(\textbf{a}, \textbf{b}) \ll 1$ can be seen by letting $p$ go to infinity in \eqref{and}, and the conditions $c_p(\textbf{a}, \textbf{a}) = 1$ and $\sigma(\textbf{b}) \leq \textbf{a} \text{ for some } \sigma \in S_n$ were shown, for instance, in \cite[Theorem 4.3]{RhSh}. \\

The  expression \eqref{and} simplifies significantly in the limit as $p \rightarrow \infty$. 
We have 
\begin{equation}\label{simplifies}
 p^{v(\textbf{a})} \omega(T_{\textbf{a}}(p ))(\textbf{x})  
= s_{\textbf{a}}(\textbf{x}) + O(1/p)
\end{equation}
where $s_{\textbf{a}}(\textbf{x})$ is the Schur polynomial \cite[Section I.3]{Mac} and the ``error term'' denotes a polynomial with coefficients of size $O(1/p)$.  In particular,
\begin{equation}\label{schur}
s_{(j, 0, \ldots, 0)}(\textbf{x}) = \sum_{\substack{\textbf{b} \in \Bbb{N}_0^n\\  |\textbf{b}| =j}}  \textbf{x}^{\textbf{b}}.
\end{equation}
 
\textbf{Interlude on partitions.} Let $\Pi(n) \subset (\Bbb{N}_0^n)^{\ast}$ be the subset of all non-increasing $n$-tuples $\textbf{a} = (a_1, \ldots, a_n)$ with $|\textbf{a}| = n$.  We order $\Pi(n)$ lexicographically.  Given  $\textbf{a}, \textbf{a}' \in \Pi(n)$, we denote by $D_{\textbf{a}', \textbf{a}}$ the number of   matrices $(c_{ij}) \in \text{Mat}(n, \Bbb{N}_0)$ with $\sum_i c_{ij} = a_j$ and $\sum_{j} c_{ij} = a'_i$, i.e., with prescribed column and row sums.  We write $D = (D_{\textbf{a}',\textbf{a}})_{\textbf{a}', \textbf{a} \in \Pi(n)}$; this is clearly a symmetric matrix. An important fact is that the matrix $D$ is invertible, and in fact has determinant $1$. This follows from \cite[I.6.7(ii)]{Mac} which says that $D$ has a Cholesky decomposition $D = A^{\top} A $ where $A = (K_{\textbf{a}, \textbf{a}'})$ is the upper uni-triangular matrix consisting of Kostka numbers $K_{\textbf{a}, \textbf{a}'}$ \cite[Section I.6]{Mac}. \\


We will be particularly interested in the Hecke operators
$$T_{[j]}( p) := T_{(j, 0, \ldots, 0)}(p )$$
with $1 \leq j \leq n$, cf.\ \eqref{am1}. We write
$$T'_{[j]}(p ) := T_{(j, \ldots, j, 0)}(p ) = \overline{T}_{[j]}( p)$$
and denote the corresponding eigenvalues by $\lambda_{[j]}(p, F)$ and $\overline{\lambda}_{[j]}(p, F)$. We start by constructing a polynomial combination of the identity. 

\begin{lemma}\label{lem1}  Let $p$ be sufficiently large.  For each $\textbf{a}= (a_1, \ldots, a_n) \in \Pi(n)$ there exists $y_{\textbf{a}} \in \Bbb{Q}$ with $  |y_{\textbf{a}}| \ll 1$ such that $$p^n \sum_{\textbf{a} \in \Pi(n)} y_{\textbf{a}} \prod_{j=1}^n T_{[a_j]}(p ) = p^{n(n+1)/2} T_{(1, \ldots, 1)}( p) = p^{n(n+1)/2}  {\rm id}.$$
\end{lemma}

\textbf{Proof.}  By Lemma \ref{basic}, \eqref{simplifies} and \eqref{schur} we have
\begin{equation}\label{formula}
 \omega(T_{[j]}( p))(\textbf{x}) = \frac{1}{p^j} \sum_{\substack{\textbf{b} \in \Bbb{N}_0^n\\ |\textbf{b}| = j}} \left( 1 + O\Bigl(\frac{1}{p}\Bigr)\right) \textbf{x}^{\textbf{b}}
 \end{equation}
for $1\leq j \leq n$. Hence for $\textbf{a} = (a_1, \ldots, a_n) \in \Pi(n)$ we have
$$p^n \omega\Bigl(\prod_{j=1}^n T_{[a_j]}(p )\Bigr)(\textbf{x}) =\prod_{j=1}^n  \omega\bigl(p^{a_j} T_{[a_j]}(p )\bigr)(\textbf{x}) =   \sum_{\textbf{a}' \in \Pi(n)} C_{\textbf{a}', \textbf{a}}\sum_{\substack{\textbf{c}\in \Bbb{N}_0^n\\ \sigma(\textbf{c}) = \textbf{a}' \text{ for some $\sigma \in S_n$}}}  \textbf{x}^{\textbf{c}}$$
where
$$C_{\textbf{a}', \textbf{a}} = D_{\textbf{a}', \textbf{a}} + O\left(\frac{1}{p}\right).$$
Indeed, if we multiply together   \eqref{formula} for various Hecke operators $T_{[a_j]}(p )$, we pick row vectors $\textbf{b}_j = (b_{j1}, \ldots, b_{jn})$ with $|\textbf{b}_j| = a_j$, and then $\textbf{x}^{\textbf{b}_1} \cdot \ldots \cdot \textbf{x}^{\textbf{b}_n} = \textbf{x}^{\textbf{c}}$ where $\textbf{c} = (c_1, \ldots, c_n)$ is given by the column sum $c_i = \sum_j b_{ji}$. 

 Recalling that $${p^{n(n+1)/2}} \omega(T_{(1, \ldots, 1)}(p ))(\textbf{x}) =   x_1 \cdot \ldots \cdot x_n , $$
we need to show that the   linear system
$$(C_{\textbf{a}', \textbf{a}})\cdot  \textbf{y} = (0, \ldots, 0, 1)^{\top} \in \Bbb{R}^{|\Pi(n)|}$$
has a (bounded) solution $\textbf{y} \in \Bbb{R}^{|\Pi(n)|}$, and this follows (for sufficiently large $p$)\footnote{One can check that the statement is true for all primes $p$} from the fact that the matrix $D = (D_{\textbf{a}', \textbf{a}})$ has determinant 1.

\begin{cor}\label{big}
If $p$ is sufficiently large, at least one of
$$\frac{ |\lambda_{[j]}(p, F)|}{p^{j(n-1)/2}} \quad (1 \leq j \leq n)$$
is $\geq (|\Pi(n)| \displaystyle \max_{\textbf{a} \in \Pi(n)} y_{\textbf{a}} )^{-1}$ and hence $ \gg 1$. 
\end{cor}
 
\textbf{Proof.} If not, then Lemma \ref{lem1} implies
$$p^{n(n+1)/2} \leq p^n \sum_{\textbf{a} \in \Pi(n)} |y_{\textbf{a}}| \prod_{j=1}^n |\lambda_{[a_j]}(p, F)| < \frac{p^n}{|\Pi(n)|} \sum_{\textbf{a} \in \Pi(n)} p^{n(n-1)/2}, $$
a contradiction.\\ 

The following lemma studies the decomposition of  $T_{[j]}(p ) T_{[j]}'(p )$. 
 
 \begin{lemma}\label{lem2}  Let $1 \leq j \leq n$. For $0 \leq i \leq j$ there exist $c_{ij} \ll 1$  such that
 $$T_{[j]}( p) T'_{[j]}(p ) = \sum_{i=0}^j c_{ij} p^{(n-1)i} T_{(2j-i, j, \ldots, j, i)}(p ) = \sum_{i=0}^j c_{ij} p^{(n-1)i} T_{(2j-2i, j-i, \ldots, j-i, 0)}(p ).$$
 \end{lemma}

\textbf{Proof.} The second equality is trivial. We prove the first. 
Taking the image under the Satake map, we obtain from \eqref{multi} an equality of polynomials
\begin{equation}\label{eq2}
  \omega\left(p^{v(j, 0, \ldots, 0) + v(j, \ldots, j, 0)} T_{[j]}(p )T'_{[j]}( p)\right)(\textbf{x}) = \sum_{\textbf{a}\in (\Bbb{N}_0^n)^{\ast}} \alpha_{\textbf{a}}\,   \omega\left(p^{v(\textbf{a})} T_{\textbf{a}}(p ) \right)(\textbf{x})
\end{equation}
for certain unique $\alpha_{\textbf{a}}\geq 0$ (they depend on $p$, but this is suppressed from the notation). Clearly the sum is restricted to $\textbf{a}$ with $|\textbf{a}| = nj$. From \eqref{multi} and the theory of invariant factors \cite[Theorem II.14]{Ne} it   follows that only those $\textbf{a} = (a_1, \ldots, a_n)$ occur in the sum for which 
\begin{equation}\label{eq1}
a_k \geq j, \quad 1 \leq k \leq n-1.
\end{equation}
We observe that  the left hand side of \eqref{eq2} is  a polynomial 
$f(x_1, \ldots, x_n)$ of degree $nj$ that is invariant under $f(x_1, \ldots, x_n) \mapsto (x_1 \cdots x_n)^{2nj}f(1/x_1, \ldots, 1/x_n)$. Hence the left hand side enjoys the same symmetry, which means
$\alpha_{(a_1, \ldots, a_n)} = \alpha_{(2j-a_n, \ldots, 2j - a_1)}$. In particular, if $\alpha_{\textbf{a}} \not= 0$, then not only \eqref{eq1} holds, but also the dual inequality $2j-a_k \geq j$ for $2 \leq k \leq n$, and hence $a_k = j$ for $2 \leq k \leq n-1$. We conclude that only $\textbf{a}$'s   of the form $(2j-i, j, \ldots, j, i)$ with $0 \leq i \leq j$ occur on the right hand side  of \eqref{eq2}. 

Now by Lemma \ref{basic}, the Satake images on both sides are polynomials with  bounded coefficients, and we can compute the numbers $\alpha_{\textbf{a}}$ by inverting the Satake map which by Lemma \ref{basic} is given by a uni-triangular matrix with bounded entries. Hence $\alpha_{\textbf{a}} \ll 1$, and so 
$$c_{ij} =  \alpha_{(2j-i, j, \ldots, j, i)}  p^{v(2j-i, j, \ldots, j, i) - v(j, 0, \ldots, 0) - v(j, \ldots, j, 0) - (n-1)i} = \alpha_{(2j-i, j, \ldots,j,  i)} \ll 1.$$


 \section{Diophantine lemmas}

We start with the following lemma \cite[Corollary 4]{BP}:

\begin{lemma}\label{lembp}
There exists a constant $A>0$ such that for each $\varepsilon, \delta, D>0$ and each quadratic polynomial $P(x,y)\in\Bbb{R}[x,y]$ whose quadratic homogeneous part is positive definite with discriminant $|\Delta|\geq D$  we have
\[
 \#\{ (x,y)\in\Bbb{Z}^2 \mid |P(x,y)|< \delta \} \ll_{D,\varepsilon} Z^\varepsilon + \delta^{1/7} Z^A
\] 
where $Z = \delta + 1 + H(P )$. 
\end{lemma}
Here $H(P )$ is the height of $P$, the size of the largest coefficient in absolute value. \\

If $1 \leq k \leq n$ and $x_1, \ldots, x_k  \in \Bbb{R}^n$, we denote by $\mathcal{V}(x_1, \ldots, x_k)$ the $k$-dimensional volume of the parallelepiped spanned by $x_1, \ldots, x_k$.

\begin{lemma}\label{constr} Let $X \geq 1$, $E > 0$, $1 \leq k \leq n$. Let $x_1, \ldots, x_k \in \Bbb{R}^n$ with entries bounded by $X$, and let $V = \mathcal{V}(x_1, \ldots, x_k)$. Let $q = (q_1, \ldots, q_k) \in \Bbb{R}^k$, and let $y \in \Bbb{R}^n$. Then the $k$ conditions $x_i^{\top} y = q_i + O(E)$ imply that after possible re-ordering of coordinates, we can decompose $y = (y_1, y_2)  \in \Bbb{R}^k \times \Bbb{R}^{n-k}$ such that $y_2 = A y_1 + b + O(F)$ with $$\| A \| \ll \frac{X^k}{V}, \quad \| b \| \ll \frac{X^{k-1}}{V} \| q \|, \quad F = \frac{X^{k-1}}{V}E.$$
\end{lemma}

\textbf{Proof.} Let $M$ be the $k \times n$-matrix whose $i$-th row is $x_i^{\top} $, $1 \leq i \leq k$. By Pythagoras, the sum of the squares of all $k \times k$-determinants of $M$ equals   $V^2$.  Hence we can find one $k\times k$-determinant  of $M$ that is $\gg V$.  
Without loss of generality let us assume that this is the leftmost determinant (otherwise exchange indices), and write $M = (M_1 \, M_2)$ where $M_1$ is a $k \times k$ matrix with determinant $\gg V$, so that by Cramer's rule  $\| M_1^{-1} \|  \ll X^{k-1}/V$. We obtain
$$y_1 = -M_1^{-1} M_2 y_2 + M_1^{-1} q + O(\| M_1^{-1}\| E)$$
which implies the lemma. 


\begin{cor}\label{count} Let $n \geq 2$. Let $Q \in \Bbb{R}^{n \times n}$ be a fixed symmetric positive definite matrix and let $X \geq 1$. Let $0 \leq k \leq n-2$, and let $x_1, \ldots, x_k \in \Bbb{Z}^n$ be linearly independent of norm $\ll X$. Let $q_0, q_1, \ldots, q_k \in \Bbb{R}$ be bounded by $X^2$ and let $0 < \delta  < X^{-N}$. Then
$$  \#\{y \in \Bbb{Z}^n \mid y^{\top}Q y = q_0 + O(X^2 \delta), x_j^{\top} Q y = q_j + O(X^2 \delta) \text{ for } 1 \leq j \leq k\} \ll X^{n-k-2 + \varepsilon}$$
if $N$ is sufficiently large. 
\end{cor}

Since $x_1, \ldots, x_k$ are integral and linearly independent, they span a parallelepiped of volume $V \geq 1$. Since $Q$ is a fixed positive definite matrix, we  have $\mathcal{V}(Qx_1, \ldots, Qx_n) \gg 1$. By Lemma \ref{constr}, the conditions $x_j^{\top} Q y = q_j + O(X^2 \delta)$ imply a decomposition $y = (y_1, y_2) \in \Bbb{R}^k \times \Bbb{R}^{n-k}$ where 
\begin{equation}\label{linear}
y_1 = Ay_2 + b + O(F), \quad \| A \| \ll X^k, \|b\| \ll X^{k+1}, F \ll X^{k+1} \delta.
\end{equation}
We insert this into the quadratic condition $y^{\top} Q y = q_0 + O(X^2\delta)$. Clearly $\| y \| \ll X$, so the error term contributes $O(X F) = O(X^{k+2} \delta)$. We write $Q$ in block notation $$Q = \left(\begin{matrix} Q_1 & Q_2\\ Q_2^{\top} & Q_3\end{matrix}\right), \quad Q_1 \in \Bbb{R}^{k \times k}, Q_2 \in \Bbb{R}^{k \times (n-k)} ,Q_3 \in \Bbb{R}^{(n-k) \times (n-k)}.$$
Inserting \eqref{linear} gives us
\begin{displaymath}
\begin{split}
y^{\top} Q y & =   (Ay_2 + b)^{\top} Q_1 (Ay_2 + b) + y_2^{\top} Q_2^{\top} (Ay_2 + b) + (Ay_2 + b)^{\top} Q_2 y_2 + y_2^{\top} Q_3 y_2 + O(X^{k+2}\delta)\\
  & = y_2^{\top}  (A^{\top} Q_1 A + Q_2^{\top} A + A^{\top} Q_2 + Q_3) y_2 + b^{\top} (Q_1A + Q_2)y + y_2^{\top} (A^{\top} Q_1 + Q_2^{\top}) b + b^{\top} Q_1 b + O(X^{k+2}\delta). 
\end{split}
\end{displaymath}
Write $Q_1 = B^{\top}B$, $Q_2 = B^{\top} C$ where $B$ is symmetric and positive definite. Write $D = BA + C$. Then, by completing the square, we have
\begin{displaymath}
\begin{split}
y^{\top} Q y  =  y_2^{\top}(D^{\top} D + Q_3 - C^{\top} C) y_2 + b^{\top} B^{\top}D y_2 + y_2^{\top} D^{\top} B b + b^{\top} Q_1 b + O(X^{k+2}\delta).
\end{split}
\end{displaymath}
All in all we see that \eqref{linear} implies
\begin{equation*}
   (Dy_2 + Bb)^{\top}(Dy_2 + Bb) +  y_2^{\top} (Q_3 - Q_2^{\top} Q_1^{-1} Q_2)y_2 = q_0 + O(X^{k+2}\delta). 
 \end{equation*}
The matrix $Q_3 - Q_2^{\top} Q_1^{-1} Q_2$ is known as the Schur complement, and it is positive definite. It is also obvious that   $(Dy_2 + Bb)^{\top}(Dy_2 + Bb) \geq 0$. Hence we can choose $n-k-2$ of the $n-k$ entries of $y_2$ freely in $X^{n-k-2}$ ways. Then we are left with an inhomogeneous binary problem whose quadratic homogeneous part $D^{\top} D + (Q_3 - Q_2^{\top} Q_1^{-1} Q_2)$ is positive definite with determinant $\gg 1$ (here we use $\det(A+B) \geq \max (\det A, \det B)$ for symmetric positive semi-definite matrices) and height $\ll X^{2k+2}$ (recall that the height is with respect to the complete polynomial, not only the homogeneous quadratic part). By Lemma \ref{lembp} we conclude that we have
 $ \ll  X^{\varepsilon}$
 choices for the remaining two entries if $N$ is sufficiently large (namely $N \geq k+2 + 14A(k+1)$ with $A$ as in Lemma \ref{lembp}).


 

\section{A general amplification scheme}\label{gen}

Let $\LL$ be a  parameter and  let $\mathcal{P}$ 
 be the set of primes $l$  in $(\LL, 2\LL]$. For $1 \leq j \leq n$ define $x_{[j]}(l) = |\lambda_{[j]}(l, F_{\varpi_0})|/\lambda_{[j]}(l, F_{\varpi_0})$  with the convention $0/0 = 0$. Recall that $F_{\varpi_0}$ is the Hecke-Maa{\ss} cusp form for which we want to prove Theorem \ref{thm1}.  For an automorphic representation $\varpi$ occurring in the spectral decomposition \eqref{spec} we define
$$A_{\varpi} :=  \sum_{j=1}^{n} \left| \sum_{l \in \mathcal{P}} \frac{x_{[j]}(l ) \lambda_{[j]}(l, F_{\varpi})}{l^{j(n-1)/2}}\right|^2 \geq 0.$$
 Then by Corollary \ref{big} we have
\begin{equation}\label{lower}
A_{\varpi_0} \gg \Bigl| \sum_{l \in \mathcal{P}} 1 \Bigr|^2 \gg |\mathcal{P}|^2. 
\end{equation}
On the other hand, by Lemma \ref{lem2}    we have
\begin{displaymath}
\begin{split}
  A_{\varpi} & =   \sum_{j=1}^{n}  \sum_{l_1\not=  l_2} \frac{x_{[j]}(l_1)\bar{x}_{[j]}(l_2) }{(l_1l_2)^{j(n-1)/2}} \lambda_{[j]}(l_1, F_{\varpi})\overline{\lambda}_{[j]}(l_2, F_{\varpi}) +  \sum_{j=1}^{n} \sum_{i=0}^j \sum_{l } \frac{|x_{[j]}(l)|^2 c_{ij}  }{l^{(n-1)(j-i)}} \lambda_{(2j-2i, j-i, \ldots, j-i, 0)}(l, F_{\varpi}).  
\end{split}
\end{displaymath}
For an integral matrix let $\Delta_j$ denote the $j$-th determinantal divisor, i.e.\ the greatest common divisor of all $j\times j$ subdeterminants. For $m, l \in \Bbb{N}$ let
$$S(m, l) := \{\gamma \in \text{Mat}(n, \Bbb{Z}) \mid \det \gamma = m,  \Delta_1 = 1, \Delta_2 = l \}.$$
Given a double coset $\Gamma \text{diag}(p^{a_1}, \ldots, p^{a_n}) \Gamma = \cup_i \Gamma M_i$ written in coset decomposition, we apply \eqref{pt} with $(x, y) = (g, 
M_i g)$. Note that all matrices in the coset decomposition of $\Gamma \text{diag}(l^{2j-2i}, l^{j-i}, \ldots, l^{j-i}, 1)\Gamma$ are in $S(l^{n(j-i)}, l^{j-i})$ and all matrices in the coset decomposition of $$\Gamma \text{diag} (l_1^j, 1, \ldots, 1) \Gamma \cdot \Gamma\text{diag}(l_2^j, \ldots, l_2^j, 1)\Gamma) =  \Gamma\text{diag}(l_1^jl_2^j, l_2^j, \ldots, l_2^j, 1)\Gamma,$$ cf.\ \cite[Prop.\ 3.2.5]{AZ}, 
are in $S(l_1^jl_2^{j(n-1)}, l_2^j)$. In this way 
we conclude from \eqref{pt} and  the definition \eqref{defHeck} of Hecke operators, as well as the lower bounds \eqref{geq1} and  \eqref{lower}   that 
\begin{equation*} 
\begin{split}
 |\mathcal{P}|^2 |F_{\varpi_0}(g)|^2 & \ll \int A_{\varpi} \tilde{f}_{\mu}(\mu_{\varpi}) |F_{\varpi}(g)|^2 d\varpi \\
 &  \ll  \sum_{j=1}^n \sum_{l_1 \not= l_2} \frac{1}{\LL^{(n-1)j}} \sum_{\gamma \in S(l_1^jl_2^{(n-1)j}, l_2^j)} |f_{\mu}(g^{-1} \gamma g)| 
 + \sum_{0 \leq i \leq j \leq n} \sum_{l} \frac{1}{\LL^{(n-1)(j-i)}} \sum_{\gamma\in S(l^{n(j-i)}, l^{j-i})} |f_{\mu}(g^{-1} \gamma g)|.
\end{split}
\end{equation*}
Recalling \eqref{spherical}, we bound the terms with $j = i$ in the last term trivially by 
$|\mathcal{P}||\textbf{c}(\mu)|^{-2}$, 
 and simplify the remaining expression as 
\begin{equation}\label{finaltrace2a}
\begin{split}
 |\mathcal{P}|^2 |F_{\varpi_0}(g)|^2 & \ll \frac{|\mathcal{P}|}{|\tilde{\textbf{c}}(\mu)|^2}+  \sum_{\nu=1}^n \sum_{l_1,  l_2 \in \mathcal{P}} \frac{1}{\LL^{(n-1)\nu}} \sum_{\gamma \in S(l_1^{\nu}l_2^{(n-1)\nu}, l_2^{\nu})} |f_{\mu}(g^{-1}  \gamma g)|. 
 \end{split}
\end{equation}
 
 Fix some large $N > 1$, and let $\delta_0  = \LL^{-N}$. We recall the notation \eqref{Cartan} for the Cartan projection and write $C_{\gamma, g} := \| C(g^{-1}\gamma g) \|$ where $\| . \|$ is some $W$-invariant norm on $\mathfrak{a}$. 
Recall from Section \ref{sec3} that $f_{\mu}$ has compact support, so only those $\gamma$ with $C_{\gamma, g} \ll 1$ contribute to the sum. For such $\gamma$ we have
$$g^{-1} \frac{\gamma}{(\det\gamma)^{1/n}} g \in K + O(C_{\gamma, g}) \quad \text{so that}   \quad  \frac{\gamma}{(\det\gamma)^{1/n}} \in g K g^{-1} + O_{\Omega}(C_{\gamma, g})$$
for $g \in \Omega$. Define $$Q = g^{-\top}g^{-1} = (q_{ij}) \in \Bbb{R}^{n \times n},$$ a fixed  positive definite symmetric matrix. Note that $gKg^{-1} = \{h \in G \mid h^{\top} Q h  = Q\}$. 
Hence a matrix $\gamma $ occurring in \eqref{finaltrace2a} with $C_{\gamma, g} \leq \delta$ satisfies 
 \begin{equation}\label{majorconditiona}
 \gamma^{\top} Q\gamma = (\det \gamma)^{2/n}Q + O((\det\gamma)^{2/n}\delta), \quad \det \gamma \asymp \LL^{n\nu}. 
 \end{equation}
Clearly all entries of $\gamma$ are bounded by $O(\LL^{\nu})$. We estimate first the contribution of the matrices $\gamma$ in \eqref{finaltrace2a} with $\delta_0 \leq C_{\gamma, g}  \ll 1$  using the bound \eqref{spherical} and the fact that the sum in \eqref{finaltrace2a} contains trivially  at most $O(\LL^{\nu n^2})$ matrices.  Hence the contribution of matrices with $C_{\gamma, g } \geq \delta_0$ is 
\begin{equation}\label{bound1a}
\ll \frac{1}{|\tilde{\textbf{c}}(\mu)|^2} ( \|\mu \| \delta_0)^{-1/2} \LL^{n^3} \ll \left(\frac{1}{|\tilde{\textbf{c}}(\mu)|^2}\right)^{1 - \frac{1}{n(n-1)}}   \LL^{ n^3+ N/2};
\end{equation}
 recall \eqref{density} for the second inequality. In particular, if $\LL$ is a sufficiently small (in terms of $n$ and $N$), but fixed power of $|\textbf{c}(\mu)|^{-2}$, this bound is acceptable. 

It remains to bound the contribution of the matrices $\gamma$ in \eqref{finaltrace2a} with $  C_{\gamma, g} \leq \delta_0$. For $0 < \delta  < 1$, $m, l \in \Bbb{N}$ let 
  \begin{equation*}\label{notation}
  \mathcal{S}(g)_{\delta}[m, l] = \{\gamma \in S(m, l) : C_{\gamma, g}  \leq \delta \}. 
  \end{equation*}
 We win if we can bound 
 \begin{equation}\label{win}
\# \mathcal{S}(g)_{\delta}[l_1^{\nu}l_2^{(n-1)\nu}, l_2^{\nu}] \ll \LL^{\nu(n-1) - \rho}
 \end{equation}
 for some fixed $\rho > 0$ and $l_1, l_2 \asymp \LL$ and $\delta = \delta_0 = \LL^{-N}$ for $N$ sufficiently large. This is the counting problem alluded to in Section \ref{sec2}.  Indeed, combining \eqref{finaltrace2a}, \eqref{bound1a} and \eqref{win}, we obtain
 \begin{equation}
  |F_{\varpi_0}(g)|^2  \ll \frac{1}{|\tilde{\textbf{c}}(\mu)|^2}\left(\frac{1}{\LL^{1-\varepsilon}} + \left(\frac{1}{|\tilde{\textbf{c}}(\mu)|^2}\right)^{ - \frac{1}{n(n-1)}}   \LL^{ n^3+ N/2} +\frac{1}{ \LL^{\rho}}\right).
  \end{equation}
  Recalling \eqref{laplace} and \eqref{density}, this clearly implies Theorem \ref{thm1}. In the final section we will show \eqref{win} in the case $n=4$, more precisely: 
  
 \begin{lemma}\label{counting1} Let $n=4$. Assume that $l = p^{\nu} \asymp \LL^{\nu}$ is a prime power with $1 \leq \nu \leq 4$, $m \asymp l^4$, and $\delta_0 = l^{-N}$ for sufficiently large $N$. Then  
$$\#\mathcal{S}(g)_{\delta_0}[m, l] \ll_g \LL^{3\nu -\frac{1}{2}+\varepsilon}.$$
\end{lemma}
  
\section{Diophantine Analysis}\label{4}

In this section we give a proof of Lemma \ref{counting1}. Throughout we specialize the formulas of the previous section to $n=4$. Let $x_j = (x_{1j}, \ldots, x_{4j})^{\top}$ denote the $j$-th column vector of $\gamma$, and let $\gamma \in \mathcal{S}(g)_{\delta}[m, l]$ with $m, l$ as in Lemma \ref{counting1}.  By \eqref{majorconditiona} we have 
\begin{equation}\label{Cij}
x_i^{\top} Q x_j =q_{ij}  m^{2/4}+ O(m^{2/4} \delta).
\end{equation}
Clearly   $x_{ij} \ll m^{1/4}$. Moreover, the determinantal divisor condition $ \Delta_2= l$ implies  in particular
\begin{equation}\label{Dij}
x_{ij}x_{i'j'} - x_{i'j}x_{ij'} \equiv 0  \, (\text{mod } l)
\end{equation}
for $1 \leq i, i', j, j' \leq 4$. 

We need to count integral vectors $x_1, \ldots, x_4 \in \Bbb{Z}^4$ satisfying the conditions \eqref{Cij} and \eqref{Dij}. We  make the   observation that the vectors $x_i$ are linearly independent (since $Q$ is non-singular) so that we can   apply Corollary \ref{count}.

Let us count the number of such matrices having one entry divisible by   $p$, say $p \mid x_{11}$ without loss of generality (otherwise exchange indices). We choose $x_{11}$ and $x_{12}$ in $(1 + m^{1/4}p^{-1}) m^{1/4}$ ways, and by Lemma \ref{lembp} we have $\ll m^{\varepsilon}$ choices for $x_{13}$ and $x_{14}$ which gives
$$\ll (1 + m^{1/4}p^{-1}) m^{1/4+\varepsilon}  \ll m^{2/4 + \varepsilon} p^{-1}$$
choices for the first column $x_1$. By Corollary \ref{count} with $k=1, 2$ and $X = m^{1/4}$ in connection with the conditions \eqref{Cij} we have $m^{1/4 + \varepsilon}$ choices for $x_2$ and $m^{\varepsilon}$ choices for $x_3$ and $x_4$, hence
 \begin{equation}\label{notcoprime}
 \ll  \LL^{3\nu-1+\varepsilon}
 \end{equation}
 choices for $\gamma$ in total. 
 
From now on we count matrices with entries coprime to $p$. Then it follows from \eqref{Dij}   that any two columns of $\gamma$ differ modulo $l$ by some multiple $a \in  (\Bbb{Z}/l\Bbb{Z})^{\ast}$. 
 Let $D := \{q_{ii} \mid 1 \leq i \leq 4\}$ and $E := \{q_{ij} \mid 1 \leq i < j \leq 4\}$. Let $$S := \{x \in \Bbb{Z}^4 \mid x^{\top} Q x \in m^{2/4}D + O(m^{2/4} \delta)\}.$$
On the set $S$ we define an equivalence relation by $x \sim y$ if $x \equiv a y$ (mod $l$) for some $a \in (\Bbb{Z}/l\Bbb{Z})^{\ast}$. We denote by $[x]$ the equivalence class of $x$ and we write
$$[x]_0 := \{y \in [x] : x^{\top} Q y \in m^{2/4}E + O(m^{2/4} \delta)\}.$$
(Note that while the set $[x]$ depends only on the equivalence class of $x$, the set $[x]_0$ depends on $x$.) Clearly, the size of $[x]$ is at most $l (1+m^{1/4}/l)^4 \asymp \LL^{\nu}$, and by applying Corollary \ref{count} with $k=0$ and $X = m^{1/4}$, the size of $S$ is $\ll m^{2/4 + \varepsilon} \asymp \LL^{2\nu+\varepsilon}$. 
 
Fix a constant $0 < \eta < 1/2$. We say that an element $x' \in S$ is \emph{$\eta$-special} if $\#[x']_0 \geq \LL^{\nu(1-\eta)}$. For $x \in S$ let $A_{\eta}(x) = \{x' \in [x] : x' \text{ is $\eta$-special}\}$. On $A_{\eta}(x)$ we define an equivalence relation by $x_1' \approx x_2'$ if $x_1'$ and $x_2'$ are linearly dependent. Each equivalence class has $O(1)$ elements. Let $A \subseteq A_{\eta}(x)/\approx$. 
Trivially,
$$  [x] \supseteq \bigcup_{x' \in A} [x']_0 .$$ For given $x'_1 \not= x_2' \in A$ we have
$$[x_1']_0 \cap [x_2']_0 \subseteq \{ y \in S \mid  (x_j')^{\top} Q y \in m^{2/4}E + O(m^{2/4}\delta), j = 1, 2\},$$
and by Corollary \ref{count} with $k=2$, the set on the right hand side has at most $O(\LL^{\varepsilon})$ elements. 
By inclusion-exclusion we get
\begin{equation*} 
\begin{split}
 \#[x] & \geq \sum_{x' \in A} \#[x']_0 - \sum_{x'_1\not= x'_2\in A} \#([x'_1]_0 \cap [x'_2]_0) \geq \#A \LL^{\nu(1-\eta)} - \#A^2 c_{\varepsilon} \LL^{\varepsilon}
\end{split}
\end{equation*}
for some constant $c_{\varepsilon}$. If  $\#(A_{\eta}(x)/\approx) \geq \#[x] \LL^{\nu(\eta-1) + \varepsilon}$, then pick any $A$ with $\#A =  \#[x] \LL^{\nu(\eta-1) + \varepsilon}$, so that
$$\#[x]  \geq \#[x] \LL^{ \varepsilon} -\#[x]^2c_{\varepsilon} \LL^{2\nu(\eta-1)+2\varepsilon}.$$
Since $\eta < 1/2$ and $\#[x] \ll \LL^{\nu}$, this is a contradiction (for $\varepsilon$ sufficiently small and $\LL$ sufficiently large), so that we conclude
$$\#A_{\eta}(x) \ll \#(A_{\eta}(x)/\approx) \leq \#[x] \LL^{\nu(\eta-1) + \varepsilon}$$
for all $x \in S$. This implies
$$\LL^{2 \nu + \varepsilon} \gg \# S = \sum_{x /\sim} \#[x]
\geq 
\LL^{\nu(1-\eta)-\varepsilon} \sum_{ x/\sim  } \#A_{\eta}(x).$$
We finally conclude
$$\sum_{\substack{x' \in S\\ x' \text{ is $\eta$-special}}} 1 = \sum_{ x/\sim  } \#A_{\eta}(x) \ll  \LL^{\nu(1+\eta) + 2\varepsilon}. $$

We are now ready for the final count of matrices $\gamma$ with entries coprime to $p$. We choose the first column in $\LL^{2\nu + \varepsilon}$ ways, at most $\LL^{\nu(1+\eta) + \varepsilon}$ of which are $\eta$-special. For the latter we choose the second column using the conditions \eqref{Cij} and Corollary \ref{count} in $\LL^{\nu + \varepsilon}$ ways, for all the others we have by definition only $\LL^{\nu(1-\eta) +\varepsilon}$ choices for the second column. Hence we have
$$\ll \LL^{3\nu-\eta +\varepsilon} + \LL^{2\nu + \eta + \varepsilon}$$
choices for the first two columns. By Corollary \ref{count} and the condition \eqref{Cij}  the other two columns are determined up to $\LL^{\varepsilon}$. Choosing $\eta = 1/2-\varepsilon$ and recalling \eqref{notcoprime}, we get
$$\#S(g)_{\delta}[m, l] \ll \LL^{3\nu +\varepsilon} (\LL^{-\nu/2 } + \LL^{-1} ), $$
which completes the proof of Lemma \ref{counting1}. 
\\

\noindent \textbf{Funding.} The first author was supported by the Volkswagen Foundation and  Starting Grant 258713 of the European Research Council.  The second author was supported by  Starting Grant 258713 of the European Research Council and  OTKA grant no. NK104183.  \\

\noindent \textbf{Acknowledgements.} The authors would like to thank Guillaume Ricotta for useful  discussions on  this paper and \cite{HRR}, and they also thank Anke Pohl and Gergely Harcos for helpful comments and encouragement. They also thank the referee for a very quick and very careful reading of the manuscript and useful suggestions.

\end{document}